\documentclass[dvips,noinfoline]{article}

\RequirePackage[lnms]{imsart}
\RequirePackage{imsartprelims}
\RequirePackage{graphicx}
\RequirePackage[colorlinks]{hyperref}
\setattribute{printed}{text}{Printed in Lithuania}

\pubyear{2008}
\volume{1}
\volumetitle{Beyond Parametrics in\\ Interdisciplinary Research:\\
Festschrift in Honor of\\ Professor Pranab K. Sen}
\firstpage{1}
\lastpage{9}
\startlocaldefs
\setattribute{journal}{name}{Collections}
\endlocaldefs

\begin{document}

\begin{titlepage}
\editor{N. Balakrishnan, Edsel A. Pe\~na and Mervyn J. Silvapulle, Editors}
\end{titlepage}

\begin{copyrightpage}
\LCCN{2007939121}
\ISBN[13]{978-0-940600-73-7}
\ISBN[10]{0-940600-73-0}
\ISSN{1939-4039}
\serieseditor{Anthony Davison}
\treasurer{Rong Chen}
\executivedirector{Elyse Gustafson}

\end{copyrightpage}

%\begin{files}
%\filelist{files.txt}
%\end{files}
%\end{document}

\makeatletter
\gdef\doi@base{http://arXiv.org/abs/}
\makeatother

\begin{contents}[doi]
\contentsline{begintocitem}{}{}
\contentsline{jobname}{imscoll1pr}{}
\contentsline{doi}{0806.4294}{}
\contentsline{title}{Preface}{v}
\contentsline{author}{N. Balakrishnan, Edsel A. Pe\~na and Mervyn J. Silvapulle}{v}
\contentsline{endtocitem}{}{}

\contentsline{begintocitem}{}{}
\contentsline{jobname}{lnms55PR}{}
\contentsline{doi}{0806.4294}{}
\contentsline{title}{Contributors to this volume}{vii}
\contentsline{author}{}{vii}
\contentsline{endtocitem}{}{}

\contentsline{begintocitem}{}{}
\contentsline{jobname}{imscoll101}{}
\contentsline{doi}{0805.2229}{}
\contentsline{title}{Pranab Kumar Sen: Life and works}{1}
\contentsline{author}{N. Balakrishnan, Edsel A. Pe\~na and Mervyn J. Silvapulle}{1}
\contentsline{endtocitem}{}{}

\contentsline{section}{Parametric inference}{}

\contentsline{begintocitem}{}{}
\contentsline{jobname}{imscoll102}{}
\contentsline{doi}{0805.2228}{}
\contentsline{title}{Analytic perturbations and systematic bias in statistical modeling and inference}{17}
\contentsline{author}{Jerzy A. Filar, Irene Hudson, Thomas Mathew and Bimal Sinha}{17}
\contentsline{endtocitem}{}{}

\contentsline{begintocitem}{}{}
\contentsline{jobname}{imscoll103}{}
\contentsline{doi}{0805.2231}{}
\contentsline{title}{Smooth estimation of mean residual life under random censoring}{35}
\contentsline{author}{Yogendra P. Chaubey and Arusharka Sen}{35}
\contentsline{endtocitem}{}{}

\contentsline{section}{Order restricted inference}{}

\contentsline{begintocitem}{}{}
\contentsline{jobname}{imscoll104}{}
\contentsline{doi}{0805.2239}{}
\contentsline{title}{Order restricted inference for comparing the cumulative incidence of a competing risk over several populations}{50}
\contentsline{author}{Hammou El Barmi, Subhash Kochar and Hari Mukerjee}{50}
\contentsline{endtocitem}{}{}

\contentsline{begintocitem}{}{}
\contentsline{jobname}{imscoll105}{}
\contentsline{doi}{0805.2242}{}
\contentsline{title}{Statistical inference under order restrictions on both rows and columns of a matrix, with an application in toxicology}{62}
\contentsline{author}{Eric Teoh, Abraham Nyska, Uri Wormser and Shyamal D. Peddada}{62}
\contentsline{endtocitem}{}{}

\contentsline{begintocitem}{}{}
\contentsline{jobname}{imscoll106}{}
\contentsline{doi}{0805.2246}{}
\contentsline{title}{On the structure of a family of probability generating functions induced by shock models}{78}
\contentsline{author}{Satrajit Roychoudhury and Manish C. Bhattacharjee}{78}
\contentsline{endtocitem}{}{}

\contentsline{section}{Bayesian inference}{}

\contentsline{begintocitem}{}{}
\contentsline{jobname}{imscoll107}{}
\contentsline{doi}{0805.2258}{}
\contentsline{title}{A Bayesian test for excess zeros in a zero-inflated power series distribution}{89}
\contentsline{author}{Archan Bhattacharya, Bertrand S. Clarke and Gauri S. Datta}{89}
\contentsline{endtocitem}{}{}

\contentsline{begintocitem}{}{}
\contentsline{jobname}{imscoll108}{}
\contentsline{doi}{0805.2264}{}
\contentsline{title}{Posterior consistency of Dirichlet mixtures of beta densities in estimating positive false discovery rates}{105}
\contentsline{author}{Subhashis Ghosal, Anindya Roy and Yongqiang Tang }{105}
\contentsline{endtocitem}{}{}

\contentsline{section}{Robust inference}{}

\contentsline{begintocitem}{}{}
\contentsline{jobname}{imscoll109}{}
\contentsline{doi}{0805.2268}{}
\contentsline{title}{Robust estimation in finite population sampling}{116}
\contentsline{author}{Malay Ghosh}{116}
\contentsline{endtocitem}{}{}

\contentsline{section}{Beyond parametrics}{}

\contentsline{begintocitem}{}{}
\contentsline{jobname}{imscoll110}{}
\contentsline{doi}{0805.2276}{}
\contentsline{title}{Estimation of population-level summaries in general semiparametric repeated measures regression models}{123}
\contentsline{author}{Arnab Maity, Tatiyana V. Apanasovich and Raymond J. Carroll}{123}
\contentsline{endtocitem}{}{}

\contentsline{begintocitem}{}{}
\contentsline{jobname}{imscoll111}{}
\contentsline{doi}{0805.2285}{}
\contentsline{title}{Smoothing-inspired lack-of-fit tests based on ranks}{138}
\contentsline{author}{Jeffrey D. Hart}{138}
\contentsline{endtocitem}{}{}

\contentsline{begintocitem}{}{}
\contentsline{jobname}{imscoll112}{}
\contentsline{doi}{0805.2292}{}
\contentsline{title}{A nonparametric control chart based on the Mann--Whitney statistic}{156}
\contentsline{author}{Subhabrata Chakraborti and Mark A. van de Wiel}{156}
\contentsline{endtocitem}{}{}

\contentsline{begintocitem}{}{}
\contentsline{jobname}{imscoll113}{}
\contentsline{doi}{0805.2300}{}
\contentsline{title}{Regression rank scores in nonlinear models}{173}
\contentsline{author}{Jana Jure\v {c}kov\'a}{173}
\contentsline{endtocitem}{}{}

\contentsline{begintocitem}{}{}
\contentsline{jobname}{imscoll114}{}
\contentsline{doi}{0805.2305}{}
\contentsline{title}{Chernoff--Savage and Hodges\textbf{--}Lehmann results for Wilks' test of multivariate independence}{184}
\contentsline{author}{Marc Hallin and Davy Paindaveine}{184}
\contentsline{endtocitem}{}{}

\contentsline{begintocitem}{}{}
\contentsline{jobname}{imscoll115}{}
\contentsline{doi}{0805.2316}{}
\contentsline{title}{\textit{\textbf{U}}-tests for variance components in one-way random effects models}{197}
\contentsline{author}{Juv\^{e}ncio S. Nobre, Julio M. Singer and Mervyn J. Silvapulle}{197}
\contentsline{endtocitem}{}{}

\contentsline{section}{Multiple comparisons}{}

\contentsline{begintocitem}{}{}
\contentsline{jobname}{imscoll116}{}
\contentsline{doi}{0805.2479}{}
\contentsline{title}{A comparison of the Benjamini\textbf{--}Hochberg procedure with some Bayesian rules for multiple testing}{211}
\contentsline{author}{Ma{\l }gorzata Bogdan, Jayanta K. Ghosh and Surya T. Tokdar }{211}
\contentsline{endtocitem}{}{}

\contentsline{begintocitem}{}{}
\contentsline{jobname}{imscoll117}{}
\contentsline{doi}{0805.2322}{}
\contentsline{title}{On the Simes inequality and its generalization}{231}
\contentsline{author}{Sanat K. Sarkar}{231}
\contentsline{endtocitem}{}{}

\contentsline{begintocitem}{}{}
\contentsline{jobname}{imscoll118}{}
\contentsline{doi}{0805.2328}{}
\contentsline{title}{Multiple testing procedures under confounding}{243}
\contentsline{author}{Debashis Ghosh}{243}
\contentsline{endtocitem}{}{}

\contentsline{section}{Mixture models}{}

\contentsline{begintocitem}{}{}
\contentsline{jobname}{imscoll119}{}
\contentsline{doi}{0805.2456}{}
\contentsline{title}{A pattern mixture model for a paired $2\times 2$ crossover design}{257}
\contentsline{author}{Laura J. Simon and Vernon M. Chinchilli}{257}
\contentsline{endtocitem}{}{}

\contentsline{begintocitem}{}{}
\contentsline{jobname}{imscoll120}{}
\contentsline{doi}{0805.2460}{}
\contentsline{title}{Projected likelihood contrasts for testing homogeneity in finite mixture models with nuisance parameters}{272}
\contentsline{author}{Debapriya Sengupta and Rahul Mazumder}{272}
\contentsline{endtocitem}{}{}

\contentsline{section}{Resampling methodology}{}

\contentsline{begintocitem}{}{}
\contentsline{jobname}{imscoll121}{}
\contentsline{doi}{0805.2470}{}
\contentsline{title}{Bootstrapping the Grenander estimator}{282}
\contentsline{author}{Michael R. Kosorok}{282}
\contentsline{endtocitem}{}{}

\contentsline{section}{Change point analysis}{}

\contentsline{begintocitem}{}{}
\contentsline{jobname}{imscoll122}{}
\contentsline{doi}{0805.2473}{}
\contentsline{title}{Ratio tests for change point detection}{293}
\contentsline{author}{Lajos Horv\'ath, Zsuzsanna Horv\'ath and Marie Hu\v skov\'a}{293}
\contentsline{endtocitem}{}{}

\contentsline{begintocitem}{}{}
\contentsline{jobname}{imscoll123}{}
\contentsline{doi}{0805.2485}{}
\contentsline{title}{On estimating the change point in generalized linear models}{305}
\contentsline{author}{Hongling Zhou and Kung-Yee Liang}{305}
\contentsline{endtocitem}{}{}

\contentsline{section}{Statistics in history}{}

\contentsline{begintocitem}{}{}
\contentsline{jobname}{imscoll124}{}
\contentsline{doi}{0805.2490}{}
\contentsline{title}{Using statistical smoothing to date medieval manuscripts}{321}
\contentsline{author}{Andrey Feuerverger, Peter Hall, Gelila Tilahun and Michael Gervers }{321}
\contentsline{endtocitem}{}{}

\contentsline{section}{Statistics in health}{}

\contentsline{begintocitem}{}{}
\contentsline{jobname}{imscoll125}{}
\contentsline{doi}{0805.2492}{}
\contentsline{title}{Sequential nonparametrics and semiparametrics: Theory, implementation and applications to clinical trials}{332}
\contentsline{author}{Tze Leung Lai and Zheng Su}{332}
\contentsline{endtocitem}{}{}

\contentsline{begintocitem}{}{}
\contentsline{jobname}{imscoll126}{}
\contentsline{doi}{0805.2496}{}
\contentsline{title}{Estimating medical costs from a transition model}{350}
\contentsline{author}{Joseph C. Gardiner, Lin Liu and Zhehui Luo}{350}
\contentsline{endtocitem}{}{}

%\contentsline{begintocitem}{}{}
%\contentsline{jobname}{imscoll127}{}
%\contentsline{doi}{10.1214/193940307000000275}{}
%\contentsline{title}{Rasch type models for measurement and analysis of health related quality of life}{364}
%\contentsline{author}{Mounir Mesbah}{364}
%\contentsline{endtocitem}{}{}

\contentsline{section}{Genetics and genomics}{}

\contentsline{begintocitem}{}{}
\contentsline{jobname}{imscoll128}{}
\contentsline{doi}{0805.2501}{}
\contentsline{title}{Correcting for selection bias via cross-validation in the classification of microarray data}{364}
\contentsline{author}{G. J. McLachlan, J. Chevelu and J. Zhu}{364}
\contentsline{endtocitem}{}{}

\contentsline{begintocitem}{}{}
\contentsline{jobname}{imscoll129}{}
\contentsline{doi}{0805.2516}{}
\contentsline{title}{An asymptotically normal test for the selective neutrality hypothesis}{377}
\contentsline{author}{Alu\'{\i }sio Pinheiro, Hildete P. Pinheiro and Samara Kiihl}{377}
\contentsline{endtocitem}{}{}

\contentsline{begintocitem}{}{}
\contentsline{jobname}{imscoll130}{}
\contentsline{doi}{0805.2523}{}
\contentsline{title}{Model selection and sensitivity analysis for sequence pattern models}{390}
\contentsline{author}{Mayetri Gupta}{390}
\contentsline{endtocitem}{}{}

\end{contents}

\newpage

\begin{preface}

  \begin{frontmatter}

    \title{Preface}

  \end{frontmatter}

  \thispagestyle{plain}

Pranab K. Sen has contributed extensively
to many areas of Statistics
   including order statistics, nonparametrics, robust inference,
   sequential methods, asymptotics, biostatistics, clinical trials,
   bioenvironmental studies and bioinformatics. His long list of over 600
   publications and 22 books and volumes along with numerous citations
   during the past 5 decades bear testimony to his work.

All three of us have had the good fortune of being associated
with him in different capacities.  He has given professional
and personal advice on many occasions to all of us, and we feel
that our lives have certainly been enriched by our association
with him.  He has been over the years a friend, philosopher and
a guide to us, and still continues to be one!

While parametric statistical inference remains ever so popular, semi-parametric,
Bayesian and nonparametric inferential methods have attracted great attention from numerous
 applied scientists because of their
weaker assumptions, which make them naturally robust and so more
appropriate in real-life applications.
This clearly signals for ``beyond
parametrics'' approaches which include nonparametrics,
semi-parametrics, Bayes methods and many others. Motivated
by this feature, and his drive in the ``beyond parametrics'' area, we
thought that it will be only appropriate for a volume in honor
of Pranab Kumar Sen to focus on this aspect of statistical
inference and its applications.  With this in mind, we have put
together this volume in order to (i) review some of the recent
developments in this direction, (ii) focus on some new
methodologies and highlight their applications, and (iii)
suggest some interesting open problems and possible new
directions for further research.

With these specific goals in mind, we invited a number of authors to contribute an article for this volume.
 These authors are not only experts in parametric, semi-parametric, Bayesian and nonparametric inferential methods, but also form a representative group
 from former students, colleagues, long-time friends, and other close professional associates of Pranab Kumar Sen.  All the articles  received have been
 properly peer reviewed according to the conditions set forth by the IMS Lecture Notes
 Editor.

It is important to mention here that this volume is not a proceedings, but rather a carefully planned volume consisting of articles that are consistent
 with the goal of highlighting developments ``beyond parametric inference'' and their applications.

Our sincere thanks to Professor Sen for having given his consent to
this venture and his advice on organisational matters whenever
we asked.
 Next, our special thanks go to all the authors who have contributed to this
 volume.
 All these authors share our respect and admiration for the various contributions and accomplishments of Pranab Kumar Sen
 and provided great cooperation during the entire course of this project.  We express our gratitude to Professors Rick Vitale  and Anthony Davison,
  the Past and Present Editors of the IMS Lecture Notes, for lending their support to this project and also for providing constant encouragement and help
  during the preparation of this volume.

We would like to thank the numerous reviewers for helping us with the
reviews of papers, Dr. Vytas Statulevi\v{c}ius for prompt assistance
related to to \LaTeX, Ms Geri Mattson for carrying out the publications
related tasks expeditiously and efficiently, and Ms Mala Raghavan and
Mr Kulan Ranasinghe for editorial assistance.
We enjoyed immensely putting this volume together, and it is with great pleasure
that we dedicate it to Pranab Kumar Sen!

\bigskip
\flushright{
\begin{tabular}{l@{}}
N. Balakrishnan\\
Hamilton, Ontario, Canada\vspace*{10pt}\\
Edesl A. Pe\~na\\
Columbia, South Carolina, USA\vspace*{10pt}\\
Mervyn J. Silvapulle\\
Melbourne, Australia\\
\end{tabular}
}

\end{preface}

\newpage
\begin{contributors}
\begin{theindex}

 % \item Anindya, R., \textit {University of Maryland Baltimore County}
  \item Apanasovich, T. V., \textit {Cornell University}

  \indexspace

  \item Balakrishnan, N., \textit{McMaster University}
  \item Bhattacharjee, M. C., \textit {New Jersey Institute of Technology}
  \item Bhattacharya, A., \textit {University of Georgia}
  \item Bogdan, M., \textit {Wroc{\l }aw University of Technology}

  \indexspace

  \item Carroll, R. J., \textit {Texas A\&M University}
  \item Chakraborti, S., \textit {University of Alabama}
  \item Chaubey, Y. P., \textit {Concordia University}
  \item Chevelu, J., \textit {University of Rennes}
  \item Chinchilli, V. M., \textit {Pennsylvania State University}
  \item Clarke, B. A., \textit {University of British Columbia}

  \indexspace

  \item Datta, G. S., \textit {University of Georgia}

  \indexspace

  \item El Barmi, H., \textit {City University of New York}

  \indexspace

  \item Feuerverger, A., \textit {University of Toronto}
  \item Filar, J. A., \textit {University of South Australia}

  \indexspace

  \item Gardiner, J. C., \textit {Michigan State University}
  \item Gervers, M., \textit {University of Toronto}
  \item Ghosal, S., \textit {North Carolina State University}
  \item Ghosh, D., \textit {Pennsylvania State University}
  \item Ghosh, J. K., \textit {Purdue University, Indian Statistical Institute}
  \item Ghosh, M., \textit {University of Florida}
  \item Gupta, M., \textit {University of North Carolina at Chapel Hill}

  \indexspace

  \item Hall, P., \textit {Australian National University}
  \item Hallin, M., \textit {Universit\' e libre de Bruxelles}
  \item Hart, J. D., \textit {Texas A\&M University}
  \item Horv\'ath, L., \textit {University of Utah}
  \item Horv\'ath, Z., \textit {University of Utah}
  \item Hu\v skov\'a, M., \textit {Charles University}
  \item Hudson, I., \textit {University of South Australia}

  \indexspace

  \item Jure\v {c}kov\'a, J., \textit {Charles University}

  \indexspace

  \item Kiihl, S., \textit {Universidade Estadual de Campinas}
  \item Kochar, S., \textit {Portland State University}
  \item Kosorok, M. R., \textit {University of North Carolina-Chapel Hill}

  \indexspace

  \item Lai, T. L., \textit {Stanford University}
  \item Liang, K.-Y., \textit {Johns Hopkins University}
  \item Liu, L., \textit {Michigan State University}
  \item Luo, Z., \textit {Michigan State University}

  \indexspace

  \item Maity, A., \textit {Texas A\&M University}
  \item Mathew, T., \textit {University of Maryland Baltimore County}
  \item Mazumder, R., \textit {Indian Statistical Institute}
  \item McLachlan, G. J., \textit {University of Queenslan}
%  \item Mesbah, M, \textit {University of Pierre et Marie Curie, Paris 6}
  \item Mukerjee, H., \textit {Wichita State University}

  \indexspace

  \item Nobre, J. S., \textit {Universidade Federal do Cear\'a}
  \item Nyska, A., \textit {Tel Aviv University}

  \indexspace

  \item Paindaveine, D., \textit {Universit\' e libre de Bruxelles}
  \item Peddada, S. D., \textit {NIEHS}
  \item Pe\~na, E., \textit {University of South Carolina}
  \item Pinheiro, A., \textit {Universidade Estadual de Campinas}
  \item Pinheiro, H. P., \textit {Universidade Estadual de Campinas}

  \indexspace

  \item Roy, A., \textit{University of Maryland Baltimore County}
  \item Roychoudhury, S., \textit {Schering Plough Research Institute}

  \indexspace

  \item Sarkar, S. K., \textit {Temple University}
  \item Sen, A., \textit{Concordia University}
  \item Sengupta, D., \textit{Indian Statistical Institute}
  \item Silvapulle, M. J., \textit{Monash University}
  \item Simon, L. J., \textit{Pennsylvania State University}
  \item Singer, J. M., \textit{Universidade de S\~ao Paulo}
  \item Sinha, B., \textit{University of Maryland Baltimore County}
  \item Su, Z., \textit{Genentech}

  \indexspace

  \item Tang, Y., \textit {SUNY Downstate Medical Center}
  \item Teoh, E., \textit {Insurance Institute for Highway Safety}
  \item Tilahun, G., \textit {University of Toronto}
  \item Tokdar, S. T., \textit {Carnegie Mellon University}

  \indexspace

  \item Van de Wiel, M. A., \textit {VU University Amsterdam}

  \indexspace

  \item Wormser, U., \textit {The Hebrew University of Jerusalem}

%  \indexspace
  \indexspace

  \item Zhou, H., \textit {US Food and Drug Administration}
  \item Zhu, J., \textit {University of Queensland}

\end{theindex}

\end{contributors}

\newpage

  \thispagestyle{plain}
\begin{figure}
  \includegraphics{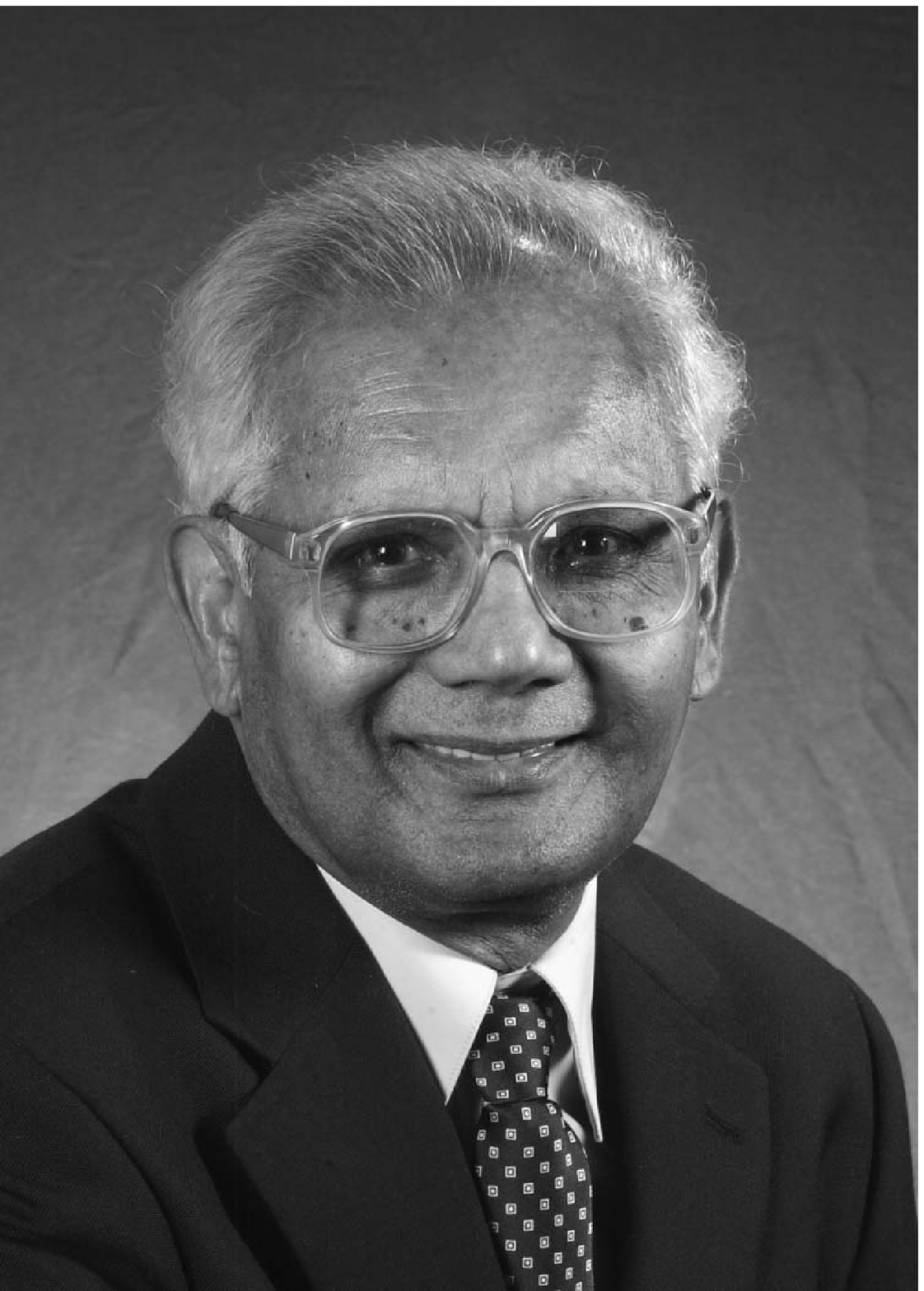}
\end{figure}

\end{document}